\documentclass[a4paper,12pt]{article}
\usepackage{amsmath, amsfonts, amssymb}
\usepackage{graphicx, color}
\usepackage[english, russian]{babel}
\makeindex

\begin{document}

\begin{center}
\textbf{BOUNDARY VALUE PROBLEM FOR A PARABOLIC-HYPERBOLIC EQUATION IN
A RECTANGULAR DOMAIN}

\textbf{Djumaklych Amanov}

Institute of Mathematics at the National University of Uzbekistan\\

e-mail: damanov@yandex.ru
\end{center}

\textbf{Abstract. }In present paper we study a boundary value problem for a mixed  parabolic-hyperbolic type equation in a rectangular domain and prove the existence of unique solution of this problem. In theory of boundary value problems for second order mixed type equations usually two conjugation conditions are in use. In this case, for mixed type equations containing hyperbolic equation in a rectangular domain for solvability of boundary value problem appears certain condition. In this paper we give three conjugation conditions.  In this case mentioned condition not appears.



\section{Introduction} 

In domain $\Omega  = \left\{ {\left( {x,t} \right):0 < x < p, - T < t < T} \right\}$ we consider equation

$$\left. \begin{array}{l}
{u_{tt}} - {u_{xx}},t > 0\\
{u_t} + {u_{xx}},t < 0
\end{array} \right\} = f\left( {x,t} \right),\eqno{(1)}$$
where $f\left( {x,t} \right)$ is given function.

We denote by  ${\Omega ^ + } = \Omega  \cap \left( {t > 0} \right),{\Omega ^ - } = \Omega  \cap \left( {t < 0} \right).$

\textbf{Problem.} Find in ${\Omega ^ + } \cup {\Omega ^ - }$ solution $u\left( {x,t} \right) \in C\left( {\overline \Omega  } \right) \cap C_{x,t}^{2,2}\left( \Omega  \right)$ of equation (1) satisfying boundary conditions
 $$u\left( {0,t} \right) = 0,u\left( {p,t} \right) = 0, - T < t < T \eqno{(2)}$$
and conjugation conditions
$$
u\left( {x, + 0} \right) = u\left( {x, - 0} \right),0 < x < p,\eqno{(3)}$$
$${u_t}\left( {x, + 0} \right) = {u_t}\left( {x, - 0} \right),0 < x < p,\eqno{(4)}$$
$${u_{tt}}\left( {x, + 0} \right) = {u_{tt}}\left( {x, - 0} \right),0 < x < p.\eqno{(5)}$$

Boundary value problems for mixed parabolic-hyperbolic type equations when hyperbolic part of domain is characteristic triangle was studied by many authors (see [1],[2],[3] and references therein). Boundary value problems for  parabolic-hyperbolic equations in rectangular domains was studied by Sabitov K.B. in [4-6] and for the higher order generalized mixed parabolic equation with fractional derivatives in [7]. The problem (1)-(5) differs from the known works of the conjugation condition (5). If we refuse condition (5) and put $u\left( {x,T} \right) = 0$ then in ${\Omega ^ + }$ we have Dirichlet type problem, the solvability of which is strongly depends on the relation of the sides of the rectangle [8]. In our case, this dependence is expressed by inequality  $\cos {\lambda _n}T + {\lambda _n}\sin {\lambda _n}T \ne 0,{\lambda _n} = \frac{{n\pi }}{p},n = 1,2,...$
the validity of which is not obvious. A similar inequality is given in [4], which proved is validity under certain restrictions on numbers $T$ and $p$ .

\section{Uniqueness of solution of problem (1)-(5)}

\textbf{Theorem 1. }Let conditions
$$\mathop {\lim }\limits_{x \to 0} x{u_x}\left( {x,t} \right) = 0,\mathop {\lim }\limits_{x \to p} \left( {p - x} \right){u_x}\left( {x,t} \right) = 0\eqno{(6)}$$
is valid, then the solution of problem (1)--(5) is unique, if it exists.

\textbf{Proof.} Let $f\left( {x,t} \right) = 0$ in $\overline \Omega  $. We will show that $u\left( {x,t} \right) = 0$ in $\overline {\Omega .} $ Following [9] consider integrals
$${\alpha _n}\left( t \right) = \int\limits_0^p {u\left( {x,t} \right){X_n}\left( x \right)dx,\ \ t > 0},\eqno{(7)}$$
$$\beta _n\left( t \right) = \int\limits_0^p u\left( x,t \right)X\left( x \right)dx,\ \ t < 0,\eqno{(8)}$$
where ${X_n}\left( x \right) = \sqrt {\frac{2}{p}} \sin {\lambda _n}x,{\lambda _n} = \frac{{n\pi }}{p},n = 1,2,...$

On the basis of (7) and (8)  we introduce the functions
$$\alpha _{n,\varepsilon }\left( t \right) = \int\limits_\varepsilon ^{p - \varepsilon } u\left( x,t \right)X_n\left( x \right)dx,\eqno{(9)}$$
$$\beta _{n,\varepsilon }\left( t \right) = \int\limits_\varepsilon ^{p - \varepsilon } u\left( x,t \right)X_n\left( x \right)dx,\eqno{(10)}$$
where $\left( {\varepsilon ,p - \varepsilon } \right) \ne \emptyset .$ Differentiating (9) two times with respect to $t$ we have
$${\alpha ''_{n,\varepsilon }}\left( t \right) = \int\limits_\varepsilon ^{p - \varepsilon } {{u_{tt}}\left( {x,t} \right){X_n}\left( x \right)dx.} $$
From homogeneous equation (1) we get ${\alpha ''_{n,\varepsilon }}\left( t \right) = \int\limits_\varepsilon ^{p - \varepsilon } {{u_{xx}}\left( {x,t} \right){X_n}\left( x \right)dx} $. Integrating by parts last integral and using conditions (6), we find
$${\alpha ''_n}\left( t \right) + \lambda _n^2{\alpha _n}\left( t \right) = 0.\eqno{(11)}$$
Analogously we get
 $${\beta '_n}\left( t \right) - \lambda _n^2{\beta _n}\left( t \right) = 0.\eqno{(12)}$$
The general solutions of equations (11), (12) have the form
$$\alpha _n\left( t \right) = a_n \cos{\lambda _n}t + b_n\sin {\lambda _n}t,$$
$$\beta _n\left( t \right) = c_n e^{\lambda _n^2t}$$.

Using conjugation conditions (3)-(5) which pass the next form
${\alpha _n}\left( 0 \right) = {\beta _n}\left( 0 \right),{\alpha '_n}\left( 0 \right) = {\beta '_n}\left( 0 \right),{\alpha ''_n}(0)= {\beta ''_n}\left( 0 \right)$ we find ${a_n} = {b_n} = {c_n} = 0.$ Consequently ${\alpha _n}\left( t \right) = {\beta _n}\left( t \right) = 0$. Then from completeness of functions ${X_n}\left( x \right)$ from (7) and (8) we obtain $u\left( {x,t} \right) \equiv 0$ in $\overline {\Omega .} $
Theorem 1 is proved.

 \section{Existence of the solution of problem (1)-(5)} 
 
 We denote $$u\left( {x,t} \right) = \left\{ \begin{array}{l}
{u^ + }\left( {x,t} \right),\ t > 0,\\
{u^ - }\left( {x,t} \right),\ t < 0.
\end{array} \right.$$
 Using Fourier's method we get the solution of problem (1)-(5) in the form
$${u^ + }\left( {x,t} \right) = \sum\limits_{n = 1}^\infty  {{X_n}} \left( x \right)\left[ {\frac{{\left( {1 - \lambda _n^2} \right){f_n}\left( 0 \right) - {{f'}_n}\left( 0 \right)}}{{\lambda _n^2\left( {\lambda _n^2 + 1} \right)}}cos{\lambda _n}t + \frac{{2{f_n}\left( 0 \right) - {{f'}_n}\left( 0 \right)}}{{{\lambda _n}\left( {\lambda _n^2 + 1} \right)}}\sin {\lambda _n}t + } \right.$$
$$\frac{1}{{{\lambda _n}}}\left. {\int\limits_0^t {{f_n}\left( \tau  \right)\sin {\lambda _n}\left( {t - \tau } \right)d\tau } } \right] \eqno{(13)}$$

$${u^ - }\left( {x,t} \right) = \sum\limits_{n = 1}^\infty  {{X_n}} \left[ {\frac{{\left( {1 - \lambda _n^2} \right){f_n}\left( 0 \right) - {{f'}_n}\left( 0 \right)}}{{\lambda _n^2\left( {\lambda _n^2 + 1} \right)}}{e^{\lambda _n^2t}} - \left. {\int\limits_t^0 {{f_n}\left( \tau  \right){e^{\lambda _n^2\left( {t - \tau } \right)}}d\tau } } \right],} \right.\eqno{(14)}$$
We consider the following derivatives
$$\begin{array}{l}
\frac{{\partial {u^ + }}}{{\partial t}}\left( {x,t} \right) = \sum\limits_{n = 1}^\infty  {{X_n}\left( x \right)\left[ {\frac{{\left( {\lambda _n^2 - 1} \right){f_n}\left( 0 \right) + {{f'}_n}\left( 0 \right)}}{{\lambda _n^2 + 1}}\sin {\lambda _n}t + } \right.} \\
+\frac{{2{f_n}\left( 0 \right) - {{f'}_n}\left( 0 \right)}}{{\lambda _n^2 + 1}}\cos {\lambda _n}t + \left. {\int\limits_0^t {{f_n}\left( \tau  \right)\cos {\lambda _n}\left( {t - \tau } \right)d\tau } } \right],
\end{array}\eqno{(15)}$$

$$\frac{{\partial {u^ - }}}{{\partial t}}(x,t) = \sum\limits_{n = 1}^\infty  {{X_n}(x)} \left[ {\frac{{2{f_n}(0) - {{f'}_n}(0)}}{{\lambda _n^2 + 1}}{e^{\lambda _n^2t}} - \int\limits_t^0 {{{f'}_n}(\tau )} {e^{\lambda _n^2(t - \tau )}}d\tau } \right],\eqno{(16)}$$

$$\frac{{{\partial ^2}{u^ + }}}{{\partial {t^2}}}\left( {x,t} \right) = \sum\limits_{n = 1}^\infty  {{X_n}\left( x \right)\left[ {\frac{{2\lambda _n^2{f_n}\left( 0 \right) + {{f'}_n}\left( 0 \right)}}{{\lambda _n^2 + 1}}\cos {\lambda _n}t + } \right.}$$
$$+\left. {\frac{{\left( {2\lambda _n^2 + 1} \right){{f'}_n}\left( 0 \right) - 2{\lambda ^2}_n{f_n}\left( 0 \right)}}{{{\lambda _n}\left( {\lambda _n^2 + 1} \right)}}\sin {\lambda _n}t + \frac{1}{{{\lambda _n}}}\int\limits_0^t {{{f''}_n}} \left( \tau  \right)\sin {\lambda _n}\left( {t - \tau } \right)d\tau } \right],\eqno{(17)}$$
$$\frac{{{\partial ^2}{u^ + }}}{{\partial {x^2}}}\left( {x,t} \right) = \sum\limits_{n = 1}^\infty  {{X_n}\left( x \right)\left[ {\frac{{2\lambda _n^2{f_n}\left( 0 \right) + {{f'}_n}\left( 0 \right)}}{{\lambda _n^2 + 1}}\cos {\lambda _n}t - {f_n}\left( t \right) + } \right.} $$
$$+\frac{{\left( {2\lambda _n^2 + 1} \right){{f'}_n}\left( 0 \right) - 2{\lambda _n}{f_n}\left( 0 \right)}}{{{\lambda _n}\left( {\lambda _n^2 + 1} \right)}}\sin {\lambda _n}t + \left. {\frac{1}{{{\lambda _n}}}\int\limits_0^t {{{f''}_n}} \left( \tau  \right)\sin {\lambda _n}\left( {t - \tau } \right)d\tau } \right],\eqno{(18)}$$
$$\frac{{{\partial ^2}{u^ - }}}{{\partial {x^2}}}\left( {x,t} \right) = \sum\limits_{n = 1}^\infty  {{X_n}\left( x \right)\left[ {\frac{{{{f'}_n}\left( 0 \right) - 2{f_n}\left( 0 \right)}}{{\lambda _n^2 + 1}}{e^{\lambda _n^2t}} + {f_n}\left( t \right) + } \right.} \left. {\int\limits_t^0 {{{f'}_n}\left( \tau  \right){e^{\lambda _n^2\left( {t - \tau } \right)}}d\tau } } \right],\eqno{(19)}$$
$$\frac{{{\partial ^2}{u^ - }}}{{\partial {t^2}}}\left( {x,t} \right) = \sum\limits_{n = 1}^\infty  {{X_n}\left( x \right)\left[ {\frac{{2\lambda _n^2{f_n}\left( 0 \right) + {{f'}_n}\left( 0 \right)}}{{\lambda _n^2 + 1}}{e^{{\lambda _n}t}} - \left. {\int\limits_t^0 {{{f''}_n}\left( \tau  \right){e^{\lambda _n^2\left( {t - \tau } \right)}}d\tau } } \right].} \right.} \eqno{20}$$

\textbf{Lemma 1.} Let  $f\left( {x,t} \right) \in C_{x,t}^{1,0}\left( {\overline \Omega  } \right),f\left( {0,t} \right) = f\left( {p,t} \right) = 0,\frac{{\partial f}}{{\partial x}} \in Li{p_\alpha }\left[ {0,p} \right],0 < \alpha  < 1$
uniformly with respect to $t$, then series
$$\sum\limits_{n = 1}^\infty  {{f_n}\left( t \right){X_n}\left( x \right)},\eqno{(21)}$$ where  $${f_n}\left( t \right) = \int\limits_0^p {f\left( {x,t} \right){X_n}\left( x \right)dx} \eqno{(22)}$$
converges absolutely and uniformly in $\overline {\Omega .} $

\textbf{Proof.} Integrating by parts integral in (22), we find
$${f_n}\left( t \right) = \frac{1}{{{\lambda _n}}}f_n^{\left( {1,0} \right)}\left( t \right)$$ where $$f_n^{\left( {1,0} \right)}\left( t \right) = \int\limits_0^p {\frac{{\partial f}}{{\partial x}}\left( {x,t} \right)\sqrt {\frac{2}{p}} \cos {\lambda _n}xdx}\eqno{(23)}$$
 according to [10]
$\left| {f_n^{\left( {1,0} \right)}\left( t \right)} \right| \le \frac{c}{{{n^\alpha }}},c = const > 0.$ Then $\left| {{f_n}\left( t \right)} \right| < \frac{{cp}}{\pi }\frac{1}{{{n^{1 + \alpha }}}}$  and series   $\sum\limits_{n = 1}^\infty  {\frac{1}{{{n^{1 + \alpha }}}}} $ converges.
Consequently, series in (21) absolutely and uniformly converges in $\overline \Omega  $ .

Lemma 1 is proved.

\textbf{Lemma 2.} If $f\left( {x,t} \right) \in C_{x,t}^{0,1}\left( {\overline \Omega  } \right),\frac{{{\partial ^2}f}}{{\partial {t^2}}} \in {L_2}\left( \Omega  \right)$ , then following estimates are valid
$$\left| {\int\limits_t^0 {f_n^{\left( k \right)}\left( \tau  \right){e^{\lambda _n^2\left( {t - \tau } \right)}}d\tau } } \right| \le \frac{{{c_k}}}{{{\lambda _n}}}{\left\| {f_n^k} \right\|_{{L_2}\left( { - T,0} \right)}},$$
$$\left| {\frac{1}{{{\lambda _n}}}\int\limits_0^t {f_n^{\left( k \right)}\left( \tau  \right)\sin {\lambda _n}\left( {t - \tau } \right)d\tau } } \right| \le \frac{{\widetilde {{c_k}}}}{{{\lambda _n}}}{\left\| {f_n^{\left( k \right)}} \right\|_{{L_2}\left( {0,T} \right)}},k = 0,1,2,$$ where ${c_k} = const > 0,\widetilde {{c_k}} = const > 0.$
The proof of this Lemma follows from H$\ddot{o}$lder inequality for integrals.

\textbf{Lemma 3.} If $f\left( {x,t} \right) \in C\left( {\overline \Omega  } \right)$ , $\frac{{\partial f}}{{\partial x}} \in {L_2}\left( \Omega  \right)$ , then series
$$\sum\limits_{n = 1}^\infty  {{X_n}\left( x \right)\frac{{2\lambda _n^2{f_n}\left( 0 \right)}}{{\lambda _n^2 + 1}}cos{\lambda _n}t} \eqno{(23)}$$
absolutely and uniformly converges in $\overline \Omega  $ .

\textbf{Proof.} Integrating by parts integral $\int\limits_0^p {f\left( {x,0} \right){X_n}\left( x \right)dx} $ we have
${f_n}\left( 0 \right) = \frac{1}{{{\lambda _n}}}f_n^{\left( {1,0} \right)}\left( 0 \right)$  and  $\frac{{2\lambda _n^2\left| {{f_n}\left( 0 \right)} \right|}}{{\lambda _n^2 + 1}} \le 2\left| {{f_n}\left( 0 \right)} \right| = \frac{2}{{{\lambda _n}}}\left| {f_n^{\left( {1,0} \right)}(0)} \right|$ .  For series (23) the series
  $\sum\limits_{n = 1}^\infty  {\frac{{\left| {f_n^{\left( {1,0} \right)}\left( 0 \right)} \right|}}{{{\lambda _n}}}} $ is majorant. Applying H$\ddot{o}$lder's inequality for sums to last series we get
$${\sum\limits_{n = 1}^\infty  {\frac{{\left| {f_n^{\left( {1,0} \right)}\left( 0 \right)} \right|}}{{{\lambda _n}}} \le {{\left( {\sum\limits_{n = 1}^\infty  {\frac{1}{{{\lambda _n}^2}}} } \right)}^{\frac{1}{2}}}{{\left( {\sum\limits_{n = 1}^\infty  {{{\left| {f_n^{\left( {1,0} \right)}\left( 0 \right)} \right|}^2}} } \right)}^{\frac{1}{2}}} \le \frac{p}{{\sqrt 6 }}\left\| {\frac{{\partial f\left( {x,0} \right)}}{{\partial x}}} \right\|} _{{L_2}\left( {0,p} \right)}}.$$
Consequently series (23) absolute and uniformly converges in ${\overline \Omega  ^ + }$ .

Lemma 3 is proved.

\textbf{Theorem 2.} Let $f\left( {x,t} \right) \in {C^1}\left( {\overline \Omega  } \right),\frac{{{\partial ^2}f}}{{\partial {t^2}}} \in {L_2}\left( \Omega  \right),f\left( {0,t} \right) = f\left( {p,t} \right) = 0,\frac{{\partial f}}{{\partial x}} \in Li{p_\alpha }[0,p]$
uniformly with respect to $t$ , then series (13)--(20) absolutely and uniformly converge in $\overline {{\Omega ^ + }} $ and $\overline {{\Omega ^ - }} $ respectively. Solutions (13) and (14) satisfy equation (1) in   ${\Omega ^ + }$ è ${\Omega ^ - }$ and conditions  (2)--(5).

\textbf{Proof.} Deducting (18) from (17) we  convinced that solutions (13) and (14) satisfied equation (1) in ${\Omega ^ + }$. Adding (16) and (19) we convinced that solutions (13) and (14) satisfied equation (1) in ${\Omega ^ - }$. From properties of functions ${X_n}\left( x \right)$ follows that  solutions (13) and (14) satisfies conditions (2). Passing to limit in (13) and (14), (15) and (16), (17) and (20) as $t \to 0,$ we convinced that solutions (13) and (14) satisfy the conditions (3)--(5).

Theorem 2 is proved.

\textbf{Remark.}  In the domain $\Omega $ if we consider mixed elliptic - hyperbolic type equation
$$\left. \begin{array}{l}
{u_{tt}} - {u_{xx}},\,\,\,t > 0\\
{u_{tt}} + {u_{xx}},\,\,\,t < 0
\end{array} \right\} = f(x,t)$$
with  conditions (2) -- (5) and $u(x,\,\, - T) = 0$, then in ${\Omega ^ + }$ the Diriclet type problem not appears.
\begin{center}
\textbf{References}
\end{center}

1. T.D. Dzhurayev,  Boundary value problems for mixed and mixed-composite type equations.Tashkent, Fan, 1979.p.240.

2. T.D. Dzhurayev, A.Sopuyev and M. Mamajanov,  Boundary value problems for equations of parabolic-hyperbolic type. Tashkent. Fan. 1986, p.220.

3. T.D.  Dzhurayev and A.Sopuyev, To the theory of partial differential equations of fourth order. Tashkent, Fan, 2000. p.144.

4. K.B. Sabitov, Tricomi problem for a mixed parabolic-hyperbolic type equation in a rectangular domain. Mat. Zametki, 2008.v.86(2). pp.273--279. (in Russian)

5. K.B. Sabitov, L.H. Rakhmanova, Initial-boundary value problem for a mixed parabolic-hyperbolic type equation in a rectangular domain. Differenc. Uravneniya
2008. v. 44, \No9, pp.1175--1181. (in Russian)

6. K.B. Sabitov, E.M.Safin. An inverse problem for a mixed parabolic-hyperbolic type equation in a rectangular domain. Izv.Vuzov.Matematika. 2010, \No4, pp.55--62. (in Russian)

7. D. Amanov and J.M. Rassias. Boundary value problems for the higher order generalized mixed-parabolic equation with fractional derivatives. Contemporary Analysis and Applied Mathematics.vol.2, \No.2, pp.198--211, 2014.

8. D. Amanov, A.Ashyralyev. Well-posedness of boundary value problems for partial differential equations of even order. Electronic Journal of Differential Equation. Vol.2014 (2014), \No 108, pp 1--8.

9. E.I.Moiseyev. About solution of a nonlocal boundary value problem by spectral method. Differenc. Uravneniya.1999. vol.35, \No8. pp. 1094-1100.(in Russian)

10.N.K. Bari. Trigonometrical  series.M.1961.(in Russian)

\end{document}